# A RENORMALIZATION FIXED POINT FOR LORENZ MAPS

BJÖRN WINCKLER

ABSTRACT. A Lorenz map is a Poincaré map for a three-dimensional Lorenz flow. We describe the theory of renormalization for Lorenz maps with a critical point and prove that a restriction of the renormalization operator acting on such maps has a hyperbolic fixed point. The proof is computer assisted and we include a detailed exposition on how to make rigorous estimates using a computer as well as the implementation of the estimates.

## 1. INTRODUCTION

The theory of renormalization has since its introduction about thirty years ago become a central tool in the study of dynamical systems. It is used, roughly speaking, to analyze maps having the property that the first-return map to some small part of the phase space resembles the original map itself. This property is usually associated with maps which lie at the "boundary of chaos", like the prototypical example of a unimodal map at the end of a period-doubling cascade. Such period-doubling cascades have been observed for maps as well as flows, but most renormalization results so far have been for one-dimensional maps (unimodal and circle maps, see e.g. [5] and [17]) with some results for higher dimensional maps (Hénon maps, see [3] and [1]). This paper contains new results for a class of so-called Lorenz flows whose dynamics can be described using the theory of renormalization.

A Lorenz flow is a three-dimensional flow possessing a singularity of saddle type with a one-dimensional unstable manifold intersecting the two-dimensional stable manifold. If the Poincaré map to a surface straddling the stable manifold can be foliated in such a way that the leaves are invariant and contracted exponentially by the Poincaré map, then the dynamics of the flow is determined by the one-dimensional map induced by the action of the Poincaré map on the leaves (see Chap. 14 of [7]). Such one-dimensional *Lorenz maps* are increasing with a jump discontinuity at the point corresponding to the stable manifold. They have been studied extensively under the additional assumption that they be expanding (see in particular [6] and [16]), but a much wider variety of dynamics is exhibited if there is also some contraction present in the form of a critical point (see [13] and [9]) and this is the situation we will consider.

The main result of this paper is that a restriction of the renormalization operator on the space of Lorenz maps has a hyperbolic fixed point, which is proved using the contraction mapping theorem on an associated operator. We use a computer to rigorously compute estimates that shows that this associated operator is indeed a contraction. This method was pioneered by Lanford [11] (see also [12]) when he proved the existence of a fixed point of the period-doubling operator on unimodal maps. However, Lanford's paper only gives a brief outline of the method he employs without an actual proof, so we have gone through quite a lot of pains to include all the missing details in this paper (many of which were borrowed from [10]).

This article is roughly divided into two halves: In Section 2 we give all the necessary definitions to state the renormalization conjecture and the main theorem, and then we go on to prove several consequences of the main theorem in Section 3. In Section 4 we describe the method used to prove the main theorem and in Section 5 we give the proof. This concludes the first half.

In the second half of the paper we give exact details on how the estimates needed to prove the main theorem are implemented on a computer.[1] The literature on this type of computer assisted proof seems to have a tradition of never including these details, most likely because it would require an order of thousands of lines of source code. We make a conscious break from this tradition and show how to implement all estimates in only 166 lines of source code.[2] The key behind this reduction in size is to use a *functional* programming language since it allows us to program in a *declarative* style: we specify *what* the program does, not *how* it is accomplished. This also has the benefit that functions cannot have *side-effects* (the output from a function only depends on its input) which makes it easier to reason about the source code. In our context this is extremely important since it means that we can check the correctness of each function in complete isolation from the rest of the source code (and a typical function is only one or two lines long which simplifies the verification of individual functions). To further minimize the risk of programming errors we choose a *strongly typed* language since these are good at catching common programming errors during compilation.

We take this opportunity to advocate the programming language Haskell for tasks similar to the one at hand — it has all the benefits mentioned above and more, but at the same time manages to produce code which runs very fast (thanks to the GHC compiler). Unfortunately, many readers will probably have had little prior exposure to Haskell and for this reason we have in Appendix E included a brief overview of Haskell as well as a table highlighting its syntax to aid the reader in understanding the source code.

## 2. Statement of the main result

In this section we state the main result, but in order to do so we first need quite a few definitions.

**Definition 2.1.** A *Lorenz map* $f$ on a closed interval $I = [l, r]$ ($l < 0 < r$) is a monotone increasing continuous function from $I \setminus \{0\}$ to $I$ such that $f(0^-) = r$, $f(0^+) = l$ (i.e. $f$ has a jump discontinuity at 0).[3]

We require that $f(x) = \varphi(|x|^\rho)$ for all $x \in (l, 0)$, where $\varphi$ is a symmetric[4] analytic map defined on some complex neighborhood of $[l, 0]$, and similarly $f(x) = \psi(|x|^\rho)$ for $x \in (0, r)$, where $\psi$ is a symmetric analytic map defined on some complex neighborhood of $[0, r]$. The maps $\varphi$ and $\psi$ are called the *analytic parts* of $f$. The constant $\rho > 1$ is called the *critical exponent* of $f$ (and is independent of $f$).

---

[1] We have structured the article so that a reader with no interest in the details of the computer estimates can skip the second half of the article.

[2] This includes: definition of main operator and its derivative (40 lines), an interval arithmetic library (30 lines), a library for computing with analytic functions (65 lines), a linear equation solver (15 lines).

[3] The notation $f(0^-)$ is shorthand for $\lim_{x \uparrow 0} f(x)$ and this limit is assumed to exist; $f(0^+)$ is defined analogously as the right-hand limit.

[4] Here 'symmetric' means $\varphi(\bar{z}) = \bar{\varphi}(z)$.

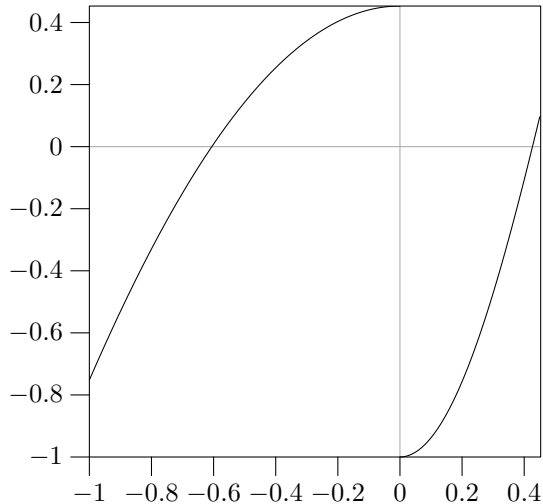

FIGURE 1. A Lorenz map (the fixed point of Theorem 2.11).

Assume $f$ is defined on $[-1, r]$ and let $g$ be a Lorenz map on $[-1, r']$ with analytic parts $(\varphi', \psi')$. We define a metric on the set of Lorenz maps by
$$d(f, g) = \max\{\|\varphi - \varphi'\|, \|\psi - \psi'\|\},$$
where $\|\cdot\|$ denotes the usual sup–norm on analytic functions. (For Lorenz maps with different domains we first perform a linear coordinate change to ensure that their domains are of the above form, then apply the above formula for the metric.)

*Remark* 2.2. The condition $\rho > 1$ ensures that $Df(x) \to 0$ as $x \to 0$ (from the left or the right) and for this reason we call $0$ the *critical point*. Because of the discontinuity at $0$ there are two *critical values*, namely $l$ and $r$.

The smoothness required in our definition of Lorenz maps is not essential for a satisfactory renormalization theory, but our current results are only in this category (which is not a big restriction since they can most likely be extended to $\mathcal{C}^r$ for $r \gtrapprox 3$ along the lines of [4] and [5]). For a discussion on what to expect when the minimum smoothness threshold is approached from below, see [2].

**Definition 2.3.** A *branch* of $f^n$ is a maximal open interval $B$ on which $f^n$ is monotone (here maximality means that if $A$ is an open interval which properly contains $B$, then $f^n$ is not monotone on $A$).

To each branch $B$ of $f^n$ we associate a word $\omega(B) = \{\sigma_0, \ldots, \sigma_{n-1}\}$ on two symbols by
$$\sigma_j = \begin{cases} 0 & \text{if } f^j(B) \subset (l, 0), \\ 1 & \text{if } f^j(B) \subset (0, r), \end{cases}$$
for $j = 0, \ldots, n - 1$.

**Definition 2.4.** A Lorenz map $f$ on $I$ is *renormalizable* if there exists a maximal interval $U \subsetneq I$ containing $0$ such that the first-return map $\tilde{f}$ to $U$ is a Lorenz map on $U$. In this situation we define the *renormalization* $\mathcal{R}f$ of $f$ as the first-return map rescaled via an increasing linear map $h : I \to U$ which takes $0$ to itself and the left endpoint of $I$ to the left endpoint of $U$:[5]

$$\mathcal{R}f = h^{-1} \circ \tilde{f} \circ h.$$

The operator $\mathcal{R}$ is called the *renormalization operator*.

*Remark* 2.5. When defined on the space of Lorenz maps with analytic branches the renormalization operator is differentiable and its derivative is a compact linear operator. This follows from the fact that $\mathcal{R}f$ only evaluates $f$ on a strict subset of the domain of $f$ (see Sections 9.4 and 9.5). On the other hand, if we only were to demand $\mathcal{C}^r$–smoothness for the branches of our Lorenz maps then $\mathcal{R}$ would no longer be differentiable (see [5] and [14, Ch. VI.1.1]).

**Definition 2.6.** Let $f$ be a renormalizable Lorenz map with associated first-return map $\tilde{f}$ and return interval $U$. Then there exists (unique) integers $a, b \geq 2$ such that $\tilde{f}|_L = f^a$ and $\tilde{f}|_R = f^b$, where $L = U \cap (l, 0)$ and $R = U \cap (0, r)$. The interval $L$ is contained in a branch $A$ of $f^a$ with associated word $\alpha = \omega(A)$, and similarly $R \subset B$ for a branch $B$ of $f^b$ with $\beta = \omega(B)$. The pair of words $(\alpha, \beta)$ is called the *type of renormalization*.

The notation $\mathcal{R}_{\alpha,\beta}$ will be used to denote the restriction of $\mathcal{R}$ to the set of Lorenz maps which have renormalizations of type $(\alpha, \beta)$.

*Remark* 2.7. In the kneading theory for unimodal maps a finite part of the kneading sequence of a point does not contain enough information to recover the exact ordering of the corresponding points on the real line. This leads to the introduction of "unimodal permutations" to describe the combinatorics of unimodal renormalization. Since Lorenz maps are strictly increasing and have no fixed points the ordering of points of a finite part of an orbit *is* completely determined by the corresponding kneading information so there is no need for this extra complication. However, we may still ask which pairs of words $(\alpha, \beta)$ give rise to valid types of renormalization. The answer to this question is given in [15]; we will not go into details as this would require more definitions, but suffice to say that there is a simple admissibility condition stated in terms of the shift operator acting on words on two symbols.

**Definition 2.8.** Let $f$ be a Lorenz map such that $\mathcal{R}^n f$ is defined for every positive integer $n$. In this situation we say that $f$ is *infinitely renormalizable*. The *combinatorial type* of $f$ is the sequence of words $\{(\alpha_0, \beta_0), (\alpha_1, \beta_1), \dots\}$, where $f$ has renormalization of type $(\alpha_0, \beta_0)$, $\mathcal{R}f$ has renormalization of type $(\alpha_1, \beta_1)$, and so on.[6] If the length of the words $\alpha_k$ and $\beta_k$ are bounded in $k$ then $f$ is said to be of *bounded combinatorial type*.

*Remark* 2.9. Let $f$ be of combinatorial type $\{\sigma_0, \sigma_1, \dots\}$, where $\sigma_k = (\alpha_k, \beta_k)$. Then $\mathcal{R}f$ has combinatorial type $\{\sigma_1, \sigma_2, \dots\}$. In other words, $\mathcal{R}$ shifts the combinatorial type to the left.

---

[5] The choice of rescaling is somewhat arbitrary (as long as it is affine) – we have chosen it so that the critical point and the left endpoint of the domain of $f$ are fixed under renormalization, whereas the right endpoint may move. Another natural choice is to fix the endpoints of the domain of $f$ but then the critical point may move.

[6] When talking about the combinatorial type we implicitly assume that it is admissible in line with the discussion of Remark 2.7.

We are now ready to state the main result, but before doing so we mention the renormalization conjecture for Lorenz maps (the statement is taken from the corresponding result for unimodal maps, see [5]):

**Conjecture 2.10** (Renormalization horseshoe)**.** *The limit set of $\mathcal{R}$ acting on the space of Lorenz maps of bounded combinatorial type is a hyperbolic Cantor set where $\mathcal{R}$ is conjugate to the full shift in a finite number of symbols, for every critical exponent $\rho > 1$.*

The above theorem represents an ultimate goal for the theory of renormalization. However, our results are much more modest in that we only prove that locally at one point in the limit set of $\mathcal{R}$ the above conjecture holds:

**Theorem 2.11** (Main Theorem)**.** *Let $\alpha = \{0,1\}$ and $\beta = \{1,0,0\}$. The restricted renormalization operator $\mathcal{R}_{\alpha,\beta}$ acting on the space of Lorenz maps with critical exponent $\rho = 2$ has a hyperbolic fixed point.*

*Proof.* This is a direct consequence of the estimates in Theorem 5.3 and the discussion in Section 4.1. □

We would like to address the question as to how far towards the renormalization conjecture our method of proof can take us. Unfortunately the answer is "not very". Theoretically, given a critical exponent $\rho > 1$ and any *periodic* combinatorial type $\{\sigma_0, \ldots, \sigma_n, \sigma_0, \ldots, \sigma_n, \ldots\}$ we could check estimates similar to those of Theorem 5.3 in order to deduce the existence of a hyperbolic fixed point. However, implementing these checks even for the simple combinatorial type at hand requires a significant effort so this does not really have any practical significance. Despite these shortcomings we still think that our current result is an important first step in the theory of renormalization of Lorenz maps.

It is also interesting to ask if any of the methods from the theory of renormalization of unimodal maps can be used to prove the renormalization conjecture for Lorenz maps. We do not know the answer to this question but it seems unlikely since the unimodal theory is based on complex analytic methods that do not have any obvious generalization to Lorenz maps (since Lorenz maps have a point of discontinuity). For this reason the renormalization theory for Lorenz maps poses new and significant difficulties. However, we can use some results from unimodal renormalization as the following remark shows.

*Remark* 2.12. The fixed point of Theorem 2.11 is the simplest non-unimodal fixed point of $\mathcal{R}$. By this we mean that if $\alpha$ and $\beta$ both have length 2 (i.e. $\alpha = \{0,1\}$, $\beta = \{1,0\}$) then the fixed point of the period-doubling operator on unimodal maps corresponds to a fixed point for $\mathcal{R}_{\alpha,\beta}$ as follows.

Let $g : [-1,1] \to [-1,1]$ be the fixed point for the period-doubling operator normalized so that $g(0) = 1$. Then $g$ is an even map that satisfies the Cvitanović–Feigenbaum functional equation
$$g(x) = -\lambda^{-1} g^2(\lambda x), \qquad \lambda = -g(1).$$

Now define a Lorenz map $f$ by $f|_{[-1,0)} = g$ and $f|_{(0,1]} = -g$. It is easy to check that the first-return map $\tilde{f}$ to $U = [-\lambda, \lambda]$ is $\tilde{f} = f^2$ and that $U$ is maximal. Thus
$$\mathcal{R}f(x) = \lambda^{-1} \tilde{f}(\lambda x) = \begin{cases} -\lambda^{-1} g^2(\lambda x) = g(x) & \text{if } x < 0, \\ \lambda^{-1} g^2(\lambda x) = -g(x) & \text{if } x > 0, \end{cases}$$

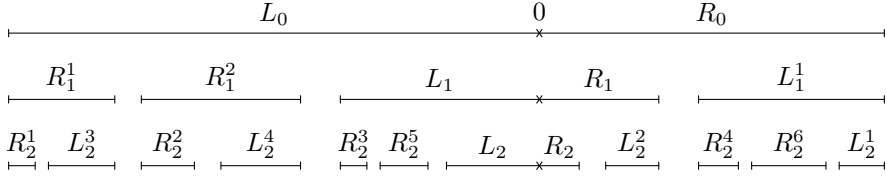

Figure 2. Illustration of the dynamical intervals of generations $0, 1, 2$ for renormalization of type $\alpha = \{0,1\}$, $\beta = \{1,0,0\}$.

which shows that $f$ is a fixed point of $\mathcal{R}_{\alpha,\beta}$.

## 3. Consequences of the main result

The existence of a hyperbolic renormalization fixed point has very strong dynamical consequences, some of which we will give a brief overview of here. Throughout this section let $\mathcal{R}$ denote the restricted renormalization operator $\mathcal{R}_{\alpha,\beta}$, where $\alpha = \{0,1\}$ and $\beta = \{1,0,0\}$, and let $f_\star$ denote the fixed point of Theorem 2.11.

**Corollary 3.1** (Stable manifold). *There exists a local stable manifold $\mathcal{W}^s_{loc}$ at $f_\star$ consisting of maps in a neighbourhood of $f_\star$ which under iteration of $\mathcal{R}$ remain in this neighbourhood and converge with an exponential rate to $f_\star$.*

*The local stable manifold extends to a global stable manifold $\mathcal{W}^s$ consisting of maps which converge to $f_\star$ under iteration of $\mathcal{R}$. If $f \in \mathcal{W}^s$ then $f$ is infinitely renormalizable.*

*Proof.* The existence of a stable manifold is a direct consequence of the stable and unstable manifold theorem. If $f$ converges to $f_\star$, then $\mathcal{R}^n f$ is defined for all $k > 0$, which is the same as saying that $f$ is infinitely renormalizable. □

We now turn to studying the dynamical properties of maps on the stable manifold. Let $f \in \mathcal{W}^s$, then the *times of closest return* $(a_n, b_n)$ are given by the recursion

$$\begin{cases} a_{n+1} = a_n + b_n, & a_1 = 2, \\ b_{n+1} = 2a_n + b_n, & b_1 = 3, \end{cases}$$

These determine the first-return interval $U_n = \text{cl}\{L_n \cup R_n\}$ for the $n$–th renormalization, by

$$L_n = \left(f^{b_n}(0^+), 0\right), \qquad R_n = \left(0, f^{a_n}(0^-)\right).$$

In other words: the first-return map $\tilde{f}_n$ to $U_n$ is given by $\tilde{f}_n(x) = f^{a_n}(x)$ if $x \in L_n$ and $\tilde{f}_n(x) = f^{b_n}(x)$ if $x \in R_n$.

Define

$$\begin{cases} L_n^k = f^k(L_n), & k = 0, \ldots, a_n - 1, \\ R_n^k = f^k(R_n), & k = 0, \ldots, b_n - 1. \end{cases}$$

The collection of these intervals (over $k$) form a pairwise disjoint collection for each $n$, called the *intervals of generation $n$* (see Figure 2).

**Theorem 3.2** (Cantor attractor). *If $f \in \mathcal{W}^s$ then the closure of the critical orbits of $f$ is a measure zero Cantor set $\Lambda_f$ which attracts almost every point in the domain of $f$.*

*Proof.* The critical orbits form the endpoints of the dynamical intervals $\{L_n^k\} \cup \{R_n^k\}$, so $\Lambda_f$ is contained in

$$\text{(1)} \qquad \bigcap_n \text{cl}\{E_n \cup F_n\}, \qquad \text{where } E_n = \bigcup_k L_n^k \text{ and } F_n = \bigcup_k R_n^k.$$

Note that $E_{n+1} \subset E_n$ and $F_{n+1} \subset F_n$.

First assume that $f = f_\star$. Then the first-return maps $\tilde{f}_n$ are all equal to $f$ itself (up to a linear change of coordinates) so the total lengths of $E_n$ and $F_n$ shrink with an exponential rate (the position of $U_{n+1}$ inside $U_n$ is the same for all $n$, so we can apply the Macroscopic and Infinitesimal Koebe principles as in the proof of the "real bounds" in [14, Ch. VI.2]). Hence the intersection (1) is a measure zero Cantor set, and consequently $\Lambda_f$ is as well.

Now, if $f$ is an arbitrary map in $\mathcal{W}^s$ then $\mathcal{R}f$ converges to $f_\star$. In other words, the first-return maps $\{\tilde{f}_n\}$ converge to $f_\star$ (up to a linear change of coordinates). Now use the same arguments as above.

For a proof that $\Lambda_f$ is an attractor, see [15]. □

**Theorem 3.3** (Rigidity)**.** *If $f, g \in \mathcal{W}^s$ then there exists a homeomorphism $h : \Lambda_f \to \Lambda_g$ conjugating $f$ and $g$ on their respective Cantor attractors. If furthermore $f, g \in \mathcal{W}_{loc}^s$, then $h$ extends to a $\mathcal{C}^{1+\alpha}$ diffeomorphism on the entire domain of $f$.*

*Proof.* Define $h(f^n(0^-)) = g^n(0^-)$ and $h(f^n(0^+)) = g^n(0^+)$. This extends continuously to a map on $\Lambda_f$ as in the proof of [14, Proposition VI.1.4].

If $f, g \in \mathcal{W}_{\text{loc}}^s$ then there exists $C > 0$ and $\lambda < 1$ such that $d(f^n, g^n) < C\lambda^n$ so we can use an argument similar to that in [14, Theorem VI.9.4] to prove the second statement. □

*Remark* 3.4 (Universality)*.* The second conclusion of Theorem 3.3 is a strong version of what is known as "metric universality": the small scale geometric structure of the Cantor attractor does not depend on the map itself (only on the combinatorial type and the critical exponent). That is, if we take two maps $f, g \in \mathcal{W}_{\text{loc}}^s$ and zoom in around the same spot on their Cantor attractors then their structures are almost identical since a differentiable map (i.e. the extended $h$) is almost linear if one zooms in closely enough.

For example, the limit of $|L_{n+1}|/|L_n|$ as $n \to \infty$ exists and is independent of $f$ (it equals the ratio $|L_2(f_\star)|/|L_1(f_\star)|$ for $f_\star$). More generally, the multifractal spectrum (and Hausdorff measure in particular) of $\Lambda_f$ does not depend on $f$ (only on $f_\star$).

We also want to mention another type of universality called "universality in the parameter plane" where the unstable eigenvalues of $DT_{f_\star}$ governs the structure of the parameter plane for families of Lorenz maps. However, in order to present the details we would need more definite information on the structure of the spectrum of $DT_{f_\star}$ so we will have to return to this discussion in another paper.

## 4. Outline of the computer assisted proof

In this section we give a brief outline of the method of proof and how to calculate rigorous estimates with a computer.

**4.1. Method of proof.** Given a Fréchet differentiable operator $T$ with compact derivative on a Banach space $X$ of analytic functions we would like to prove that $T$ has a hyperbolic fixed point. The main tool is the following consequence of the contraction mapping theorem:

**Proposition 4.1.** *Let $\Phi$ be a Fréchet differentiable operator on a Banach space $X$, let $f_0 \in X$, and let $B_r(f_0) \subset X$ be the closed ball of radius $r$ centered on $f_0$. If there are positive numbers $\varepsilon, \theta$ such that*

(1) $\|D\Phi_f\| < \theta$, for all $f \in B_r(f_0)$,
(2) $\|\Phi f_0 - f_0\| < \varepsilon$,
(3) $\varepsilon < (1-\theta)r$,

*then there exists $f_\star \in B_r(f_0)$ such that $\Phi f_\star = f_\star$ and $\Phi$ has no other fixed points inside $B_r(f_0)$. Furthermore, $\|f_\star - f_0\| < \varepsilon/(1-\theta)$.*

Our strategy is to find a good approximation $f_0$ of a fixed point of $T$ and then use a computer to verify that the conditions on $r, \varepsilon, \theta$ hold if $r$ is chosen small enough. Unfortunately, this is not possible for $T$ itself since in our case it is not a contraction (it has expanding eigenvalues) so first we have to turn $T$ into a contraction without changing the set of fixed points. This is done by using Newton's method to solve the equation $Tf - f = 0$, which results in the iteration

$$f \mapsto f - (DT_f - I)^{-1}(Tf - f),$$

where $I$ denotes the identity operator on $X$. The operator we use is a slight simplification of this, namely

$$\Phi f = f - (\Gamma - I)^{-1}(Tf - f) = (\Gamma - I)^{-1}(\Gamma - T)f,$$

where $\Gamma$ is a finite-rank linear approximation of $DT_{f_0}$ (chosen so that $\Gamma - I$ is invertible). The operator $\Phi$ *is* a contraction if $f_0$ and $\Gamma$ are chosen carefully.[7] Note that $\Phi f = f$ if and only if $Tf = f$, so once we verify that the conditions of Proposition 4.1 hold for $\Phi$ it follows that $T$ has a fixed point.

To prove hyperbolicity we need to do some extra work. The derivative of $\Phi$ is

$$D\Phi_f = (\Gamma - I)^{-1}(\Gamma - DT_f).$$

At this stage we will have already checked that the norm of this is bounded from above by 1. By strengthening this estimate to

$$\|(\Gamma - e^{it}I)^{-1}(\Gamma - DT_f)\| < 1, \qquad \forall t \in \mathbb{R}, \ \forall f \in B_r(f_0),$$

we also get that $DT_{f_\star}$ is hyperbolic at the fixed point $f_\star$. To see this, assume that $e^{it}$ is an eigenvalue of $DT_{f_\star}$ with eigenvector $h$ normalized so that $\|h\| = 1$. Then

$$\|(\Gamma - e^{it}I)^{-1}(\Gamma - DT_{f_\star})h\| = \|(\Gamma - e^{it}I)^{-1}(\Gamma - e^{it}I)h\| = \|h\| = 1,$$

which is impossible. Since $DT$ was assumed to be compact we know that the spectrum is discrete, so the lack of eigenvalues on the unit circle implies hyperbolicity.

4.2. **Rigorous computer estimates.** In order to verify the above estimates on a computer we are faced with two fundamental problems: (i) arithmetic operations on real numbers are carried out with finite precision which leads to rounding problems, (ii) the space of analytic functions is infinite dimensional so any representation of an analytic function needs to be truncated.

The general idea to deal with these problems is to compute with sets which are guaranteed to contain the exact result instead of computing with points: real numbers are replaced with intervals, analytic functions are replaced with rectangle sets $A_0 \times \cdots \times A_k \times \{C\}$ in $\mathbb{R}^n$ representing all functions of the form

$$\{a_0 + \cdots + a_k z^k + z^d h(z) \mid a_j \in A_j, \ j = 0, \ldots, k, \ \|h\| \leq C\},$$

where $\{A_j\}$ are intervals. This takes care of the truncation problem and the rounding problem is taken care of roughly by "rounding outwards" (lower bounds are rounded down, upper bounds are rounded up). Once these set representations have been chosen we lift operations on points to operations on sets. Since the form of these sets are most likely not preserved by such operations, this lifting involves finding bounds by sets of the chosen form (e.g. if $F$ and $G$ are rectangle sets of analytic functions and we want to lift composition of functions, then we have to find a rectangle set which contains the set $\{f \circ g \mid f \in F, \ g \in G\}$.)

---

[7] Here is how to choose $f_0$ and $\Gamma$: use the Newton iteration on polynomials of some fixed degree to determine $f_0$ and set $\Gamma = DT_{f_0}$. The hardest part is finding an initial guess such that the iteration converges.

Section 7 contains all the details for computing with intervals and Section 9 contains all the details for computing with rectangle sets of analytic functions.

Let us make one final remark concerning the evaluation of the operator norm of a linear operator $L$ on the space of analytic functions. In order to get good enough bounds on the estimate of the operator norm we will use the $\ell^1$–norm on the Taylor coefficients of analytic functions. The reason for this is that estimating the operator norm with

$$\|L\| = \sup_{\|f\| \leq 1} \|Tf\|$$

will usually result in really bad estimates. With the $\ell^1$–norm, if we think of $L$ as an infinite matrix (in the basis $\{z^k\}$), the operator norm is found by taking the supremum over the norms of the columns of this matrix, that is

$$\|L\| = \sup_{k \geq 0} \|L\xi_k\|, \qquad \xi_k(z) = z^k.$$

Evaluating the norms of columns gives much better estimates and for this reason we choose this norm. See Section 9.11 for the specifics.

## 5. Proof of the main theorem

First we restate the definition of the restricted renormalization operator, then we change coordinates and restate the main theorem.

### 5.1. Definition of the operator.
From now on we fix the domain of our Lorenz maps to some interval $[-1, r]$. The right endpoint cannot be fixed since it generally changes under renormalization (we will soon change coordinates so that the domain is fixed).

Instead of dealing with functions with a discontinuity we represent a Lorenz map $F$ by a pair $(f, g)$, with $f : [-1, 0] \to [-1, r]$, $f(0) = r$, and $g : [0, r] \to [-1, r]$, $g(0) = -1$.

With this notation, the first-return map to some interval $U$ will be of the type $(F^a, F^b)|_U$ if $F$ is renormalizable. For the type $\alpha = \{0, 1\}$, $\beta = \{1, 0, 0\}$, we can be more precise: in this case $a = 2$, $b = 3$ and the first-return map is of the form $(g \circ f, f \circ f \circ g)|_U$ if it is renormalizable.

Let $T$ denote the restricted renormalization operator $\mathcal{R}_{\alpha,\beta}$, and fix the critical exponent $\rho = 2$. If $T(f,g) = (\hat{f}, \hat{g})$ then $T$ is defined by

$$\hat{f}(z) = \lambda^{-1} g \circ f(\lambda z),$$
$$\hat{g}(z) = \lambda^{-1} f \circ f \circ g(\lambda z),$$
$$\lambda = -f^2(-1).$$

5.2. **Changing coordinates.** To ensure the correct normalization ($g(0) = -1$) and the correct critical exponent ($\rho = 2$) we make two coordinate changes and calculate how the operator $T$ transforms. We will also carefully choose the domain of $T$ so that all compositions are well-defined (e.g. $\lambda z$ is in the domain of $f$ etc.). This is checked automatically by the computer (and also shows that $T$ is differentiable with compact derivative, since $f$ and $g$ are analytic). Finally, it is important to realize that the choice of coordinates may greatly affect the operator norm of the derivative; not every choice will give a good enough estimate.

The domain of $T$ is chosen to be contained in the set of Lorenz maps $(f,g)$ with representation $f(z) = \phi(z^2)$ and $g(z) = \psi(z^2)$, where $\phi$ and $\psi$ have domains $\{z : |z-1| < s\}$ and $\{z : |z| < t\}$, respectively (the constants $s$ and $t$ will soon be specified). Rewriting $T$ in terms of $\phi$ and $\psi$ gives

$$\hat{\phi}(z) = \lambda^{-1} \psi(\phi(\lambda^2 z)^2)$$
$$\hat{\psi}(z) = \lambda^{-1} \phi(\phi(\psi(\lambda^2 z)^2)^2)$$
$$\lambda = -\phi(\phi(1)^2)$$

This coordinate change ensures the correct critical exponent.

The next coordinate change is to fix the normalization and also to bring the domain of all functions to the unit disk. Fixing the normalization has the benefit that the error involved in the evaluation of $\lambda$ is minimized (since we only need to evaluate $f$ close to $z = 0$, see Section 9.8). Changing all domains to the unit disk simplifies the implementation of the computer estimates.

**Definition 5.1.** Define $X$ to be the Banach space of symmetric (with respect to the real axis) analytic maps on the unit disk with finite $\ell^1$–norm. That is, if $f \in X$ then $f(z) = \sum a_k z^k$ with $a_k \in \mathbb{R}$ and $\|f\| = \sum |a_k| < \infty$.

**Definition 5.2.** Define $Y = X \times X$ with the norm $\|(f,g)\|_Y = \|f\|_X + \|g\|_X$ and with linear structure defined by $\alpha(f,g) + \beta(f',g') = (\alpha f + \beta f', \alpha g + \beta g')$. Clearly $Y$ is a Banach space (since $X$ is).

Change coordinates from $\phi, \psi$ to $(f,g) \in Y$ (note that $f$ and $g$ are not the same as above) as follows

$$\phi(z) = f([z-1]/s),$$
$$\psi(z) = -1 + z \cdot g(z/t),$$

where we will choose $s = 2.2$ and $t = 0.5$. Rewriting $T$ in terms of $f$ and $g$ gives

$$\hat{f}(w) = \lambda^{-1} \left\{ -1 + f\left(\lambda^2 \left[w + \tfrac{1}{s}\right] - \tfrac{1}{s}\right)^2 \cdot g\left(\tfrac{1}{t} \cdot f\left(\lambda^2 \left[w + \tfrac{1}{s}\right] - \tfrac{1}{s}\right)^2\right) \right\},$$
$$\hat{g}(w) = \frac{1}{tw} \left\{ 1 + \lambda^{-1} f\left(\tfrac{1}{s} \cdot f\left(\lambda^2 twg\left(\lambda^2 w\right) \cdot \left[\lambda^2 twg\left(\lambda^2 w\right) - 2\right] \cdot \tfrac{1}{s}\right)^2 - \tfrac{1}{s}\right) \right\},$$
$$\lambda = -f\left([f(0)^2 - 1]/s\right).$$

This is the final form of the operator that will be studied.

5.3. **Computing the derivative.** In order to simplify the computation of the derivative of $T$ we break the computation of $T$ down into several steps as follows:

$$\begin{aligned}
p_f(w) &= \lambda^2 \cdot (w + s^{-1}) - s^{-1} & p_g(w) &= \lambda^2 w \\
f_1 &= f \circ p_f & g_1 &= g \circ p_g \\
f_2 &= f_1^2 & g_2 &= t \cdot p_g \cdot g_1 \\
f_3 &= f_2/t & g_3 &= g_2 \cdot (g_2 - 2)/s \\
f_4 &= g \circ f_3 & g_4 &= f \circ g_3 \\
f_5 &= -1 + f_2 \cdot f_4 & g_5 &= (g_4^2 - 1)/s \\
f_6 &= f_5/\lambda & g_6 &= f \circ g_5 \\
& & g_7 &= g_6/\lambda \\
& & g_8(w) &= (g_7(w) + 1)/(t \cdot w)
\end{aligned}$$

With this notation we have that $T(f, g) = (f_6, g_8)$. Note that the result of $g_7(w)+1$ is a function with zero as constant coefficient so in the implementation of $g_8$ we will not actually divide by $w$, instead we will 'shift' the coefficients to the left.

It is now fairly easy to derive expressions for the derivative. If $f$ is perturbed by $\delta f$ and $g$ is perturbed by $\delta g$, then the above functions are perturbed as follows:

$$\begin{aligned}
\delta p_f(w) &= 2 \cdot \lambda \cdot \delta\lambda \cdot (w + s^{-1}) & \delta p_g(w) &= 2 \cdot \lambda \cdot \delta\lambda \cdot w \\
\delta f_1 &= Df \circ p_f \cdot \delta p_f + \delta f \circ p_f & \delta g_1 &= Dg \circ p_g \cdot \delta p_g + \delta g \circ p_g \\
\delta f_2 &= 2 f_1 \delta f_1 & \delta g_2 &= t \cdot (\delta p_g \cdot g_1 + p_g \cdot \delta g_1) \\
\delta f_3 &= \delta f_2/t & \delta g_3 &= \delta g_2 \cdot (g_2 - 2)/s + g_2 \cdot \delta g_2/s \\
\delta f_4 &= Dg \circ f_3 \cdot \delta f_3 + \delta g \circ f_3 & \delta g_4 &= Df \circ g_3 \cdot \delta g_3 + \delta f \circ g_3 \\
\delta f_5 &= \delta f_2 \cdot f_4 + f_2 \cdot \delta f_4 & \delta g_5 &= 2 \cdot g_4 \cdot \delta g_4/s \\
\delta f_6 &= \delta f_5/\lambda - f_5 \cdot \delta\lambda/\lambda^2 & \delta g_6 &= Df \circ g_5 \cdot \delta g_5 + \delta f \circ g_5 \\
& & \delta g_7 &= \delta g_6/\lambda - g_6 \cdot \delta\lambda/\lambda^2 \\
& & \delta g_8(w) &= \delta g_7(w)/(t \cdot w)
\end{aligned}$$

With this notation we have that $DT_{(f,g)}(\delta f, \delta g) = (\delta f_6, \delta g_8)$.

5.4. **New statement of the main theorem.** We now state the main theorem in the form it will be proved. The discussion in Section 4.1 shows how this result can be used to deduce Theorem 2.11.

**Theorem 5.3.** *There exists a Lorenz map $F_0$ and a matrix $\Gamma$ such that the simplified Newton operator $\Phi = (\Gamma - I)^{-1}(\Gamma - T)$ is well-defined and satisfies:*
  (1) $\|D\Phi_F\| < 0.2$, *for all* $\|F - F_0\| \leq 10^{-7}$,
  (2) $\|\Phi F_0 - F_0\| < 5 \cdot 10^{-9}$.
  (3) $\|(\Gamma - e^{it}I)^{-1}(\Gamma - DT_F)\| < 0.9$, *for all* $t \in \mathbb{R}$, $\|F - F_0\| \leq 10^{-7}$.

*Proof.* The remainder of this article is dedicated to rigorously checking the first two estimates with a computer. The third estimate is verified by covering the unit circle with small rectangles and using the same techniques as in the first two estimates to get rigorous upper bounds on the operator norm. However, we have left out the source code for this estimate to keep the page count down and also because the running time of the program went from a few seconds to several hours (we had to cover the circle with 50000 rectangles in order for the estimate to work). □

*Remark* 5.4. The approximate fixed point $F_0$ and approximate derivative $\Gamma$ at the fixed point are found by performing a Newton iteration eight times on an initial guess (which was found by trial-and-error). We will not spend too much time talking about these approximations but they could potentially be used to compute e.g. the Hausdorff dimension of the Cantor attractor of maps on the local stable manifold.

We did however compute the eigenvalues of $\Gamma$ and it turns out that $\Gamma$ has two simple expanding eigenvalues $\lambda_s \approx 23.36530$ and $\lambda_w \approx 12.11202$, and the rest of the spectrum is strictly contained in the unit disk. Since $\Gamma$ is a good approximation of $DT_{f_\star}$ and both operators are compact it seems clear that the spectrum of $DT_{f_\star}$ also must have exactly two unstable eigenvalues.

Lanford [12] claims that in the case of the period-doubling operator if an analog of the third estimate of Theorem 5.3 holds and "if $\Gamma$ has spectrum inside the unit disk except for a single simple expanding eigenvalue, then the same will be true for $DT_{f_\star}$." It seems plausible that a similar statement holds in the present situation with two simple expanding eigenvalues but have not yet managed to prove this (it is easy to see that if $\Gamma$ and $DT_{f_\star}$ were both diagonal then the third estimate would imply that they have the same number of unstable eigenvalues).

## 6. Implementation of computer estimates

In this section we implement the main operator and compute the estimates of Theorem 5.3. Before reading this section it may be a good idea to take a quick glance at the beginning of Section 9 in order to understand the way analytic functions are represented. It may also be helpful to use Table 1 in Appendix E to look up unfamiliar syntax in the source code.

6.1. **The main program.** To begin with we import two functions from the standard library that will be needed later:

```
1  import Data.List (maximumBy,transpose)
```

The entry point of the program is the function `main`, all that is done here is to print the result of the computations to follow:

```
2  main = do putStrLn $ "radius     = " ++ show beta
            putStrLn $ "|Phi(f)-f| < " ++ show eps
            putStrLn $ "|DPhi|     < " ++ show theta
```

The initial guess[8] is first improved by iterating a polynomial approximation[9] of the operator $\Phi$ eight times (the derivative is recomputed in each iteration, so this is a Newton iteration):

```
5  approxFixedPt = iterate (\t -> approx $ opPhi (gamma t) t) guess !! 8
```

Compute the approximation $\Gamma$ of the derivative $DT$ at the approximate fixed point:

```
6  approxDeriv = gamma approxFixedPt
```

---

[8]See Appendix C
[9]See Section 9.10

Compute an upper bound on the distance[10] between the approximate fixed point and its image under $\Phi$:

```
7   eps = upper $ dist approxFixedPt (opPhi approxDeriv approxFixedPt)
```

Construct a ball[11] of radius $\beta$ centered around the approximate fixed point and then compute the supremum of the operator norm[12] of the derivative on this ball:

```
8   theta = opnorm $ opDPhi approxDeriv (ball beta approxFixedPt)
```

The rest of this section will detail the implementation of the operator $\Phi$ and its derivative. The generic routines for rigorous computation with floating point numbers and analytic functions are discussed in the sections that follow. All input to the program (`d`, `sf`, `sg`, `guess`, `beta`) is collected in Appendix C. Instructions on how to run the program and the output it produces is given in Appendix D.

6.2. **The main operator.** The operator $T$ is computed in a function called `mainOp` which takes a Lorenz map $(f,g) \in Y$ and a sequence of tangent vectors $\{(\delta f_k, \delta g_k) \in Y\}_{k=1}^n$ and returns $(T(f,g), \{DT_{(f,g)}(\delta f_k, \delta g_k)\}_{k=1}^n)$. We perform both computations in one function since the derivative uses a lot of intermediate results from the computation of $T(f,g)$.

Given a Lorenz map (`f`,`g`) and a list of tangent vectors `ds`, first compute $f_6$ and $g_8$ as in Section 5.3 and split the result so that the polynomial parts have degree at most $d-1$. Then compute the derivatives and return the result of these two computations in a pair:

```
9   mainOp (f,g) ds = ((split d f6,split d g8), mainDer ds) where
        l  = lambda f
        pf = F [(l^2-1)/sf,l^2] 0 ; g1 = compose g pg
        pg = F [0,l^2] 0          ; g2 = pg * g1 .* sg
        f1 = compose f pf         ; g3 = g2 * (g2 - 2) ./ sf
        f2 = f1^2                 ; g4 = compose f g3
        f3 = f2 ./ sg             ; g5 = (g4^2 - 1) ./ sf
        f4 = compose g f3         ; g6 = compose f g5
        f5 = -1 + f2*f4           ; g7 = g6 ./ l
        f6 = f5 ./ l              ; g8 = lshift g7 ./ sg
```

The actual computation of the derivative is performed next inside a local function to `mainOp`. If there are no tangent vectors, no computation is performed:

```
19      mainDer [] = []
```

Otherwise, recurse over the list of tangent vectors and compute $\delta f_6$ and $\delta g_8$ and again split the result so that the polynomial parts have degree at most $d-1$:

```
20      mainDer ((df,dg):ds) = (split d df6,split d dg8) : mainDer ds where
            dl  = dlambda f df
            dpf = F ([2*l*dl]*[1/sf,1]) 0 ; dg1 = dcompose g pg dg dpg
            dpg = F [0,2*l*dl] 0          ; dg2 = (dpg*g1 + pg*dg1) .* sg
            df1 = dcompose f pf df dpf    ; dg3 = 2*(dg2*g2 - dg2) ./ sf
            df2 = 2*f1*df1                ; dg4 = dcompose f g3 df dg3
            df3 = df2 ./ sg               ; dg5 = 2*g4*dg4 ./ sf
```

---

[10]See Section 9.3
[11]See Section 9.12
[12]See Section 9.11

```
            df4 = dcompose g f3 dg df3    ; dg6 = dcompose f g5 df dg5
            df5 = df2*f4 + f2*df4         ; dg7 = dg6./l - g6.*(dl/l^2)
            df6 = df5./l - f5.*(dl/l^2)   ; dg8 = lshift dg7 ./ sg
```

Note that the constants $s$ and $t$ of Section 5.3 are called `sf` and `sg` respectively in the source code.

The above function can be used to compute the action of $T$ by passing an empty list of tangent vectors and extracting the first element of the returned pair:

30  `opT fg = fst $ mainOp fg []`

Similarly, we can evaluate $DT$ by extracting the second element:

31  `opDT fg ds = snd $ mainOp fg ds`

Using this function we compute an approximation $\Gamma$ of $DT_{(f,g)}$ by evaluating the derivative at the $2d$ first basis vectors[13] of $Y$ and approximating the result with polynomials and packing them into a $2d \times 2d$ matrix (transposing the resulting matrix is necessary because the linear algebra routines[14] we use require the matrix to be stored in row-major order):

32  ```
gamma fg = transpose $ map (interleavePoly . approx)
                    $ opDT fg (take (2*d) basis)
```

Finally, the operator $\Phi$ (and its derivative) is implemented by taking a Newton step[15] with $T$ (for convenience we pass the approximate derivative as the parameter `m`):

34  ```
opPhi  m x    = newton m (opT x) x
opDPhi m x ds = [ newton m a b | (a,b) <- zip (opDT x ds) ds ]
```

6.3. **The rescaling factor.** With our choice of coordinates the rescaling factor $\lambda$ only depends on $f$ (and not on $g$):

$$\lambda(f) = -f\big([f(0)^2 - 1]/s\big)$$

The implementation is straightforward:

36  `lambda f = -eval f (((eval f 0)^2-1)/sf)`

If $f \in X$ is perturbed by $\delta f \in X$ then $\lambda$ is perturbed by $\delta\lambda$, where

$$\delta\lambda = -2 \cdot s^{-1} \cdot f(0) \cdot \delta f(0) \cdot Df\big([f(0)^2 - 1]/s\big) - \delta f\big([f(0)^2 - 1]/s\big).$$

Derivative evaluation has to be handled carefully since we are using the $\ell^1$–norm, see Section 9.4. If $y = [f(0)^2 - 1]/s$, then $y$ lies in the closed disk of radius $|y|$ but since we need to evaluate the derivative on an open disk we first enlarge the bound on $|y|$ to get the radius $\mu$ and then evaluate the derivative on this slightly larger disk:

37  ```
dlambda f df = -2/sf * f0 * eval df 0 * eval (deriv mu f) y - eval df y
    where f0 = eval f 0
          y  = (f0^2 - 1)/sf
          mu = enlarge $ abs y
```

---

[13] See Section 9.11.
[14] See Appendix A.
[15] See Section 9.13.



## 7. Computation with floating point numbers

We discuss how to control rounding and avoid overflow and underflow when computing with floating point numbers. We show how to lift operations on floating point numbers to intervals and then how to bound these operations.

7.1. **Safe numbers.** In order to avoid overflow and underflow during the course of the proof we restrict all computations to the set of *safe numbers* (see [10]) which we define as the subset of double precision floating point numbers (referred to as *floats* from now on) $x$ such that $x = 0$ or $2^{-500} < |x| < 2^{500}$. We say that $y$ is a *safe upper bound* on $x \neq 0$ iff $x < y$ (strict inequality) and $y$ is a safe number; safe lower bounds are defined analogously. If $x = 0$, then $y = 0$ is both a safe upper and lower bound and there are no other safe bounds on $x$ (this will make sense after reading the assumption below).

Safe numbers allow us to perform rigorous computations on any computer conforming to the IEEE 754 standard since such a computer must satisfy the following assumption:[16]

**Assumption.** *Let $\bar{x}$ be a float resulting from an arithmetic operation on safe numbers performed by the computer and let $x$ be the exact result of the same operation. If $\bar{x} \neq 0$ then either $\bar{x} = x^-$ or $\bar{x} = x^+$, where $x^-$ is the largest float such that $x^- \leq x$ and $x^+$ is the smallest float such that $x \leq x^+$. Furthermore, $\bar{x} = 0$ if and only if $x = 0$.*

Under this assumption we know that the exact result must lie within any safe upper and lower bounds on $\bar{x}$, and we know that when the computer returns a result of 0 then the computation must be exact.

Given a float $x$ we now show how to find safe upper and lower bounds on $x$.

Check if a number is safe:

```
41  isSafe x = let ax = abs x in x == 0 || (ax > 2^^(-500) && ax < 2^500)
```

Use this function to assert that a number is safe, abort the program otherwise:

```
42  assertSafe x | isSafe x  = x
                 | otherwise = error "assertSafe: not a safe number"
```

Given a float we can 'step' to an adjacent float as follows:

```
44  stepFloat n 0 = 0
    stepFloat n x = let (s,e) = decodeFloat x in encodeFloat (s+n) e
```

That is, `stepFloat 1 x` is the smallest float $y$ larger than $x$ and `stepFloat (-1) x` is the largest float smaller than $x$, unless $x = 0$ in which case $x$ is returned. (The function `decodeFloat` converts a float to the form $s \cdot 2^e$, where $s, e \in \mathbb{Z}$, and `encodeFloat` converts back to a float.)

Now finding a safe upper or lower bound is easy, just step to the next float and assert that it is safe:

```
46  safeUpperBound = assertSafe . stepFloat 1
    safeLowerBound = assertSafe . stepFloat (-1)
```

---

[16]This statement follows from: (1) the fact that IEEE 754 guarantees correct rounding, (2) the result of an arithmetic operation on safe numbers is a normalized float so silent underflow to zero cannot occur.



7.2. **The Scalar data type.** The `Scalar` data type represents safe lower and upper bounds on a number:

```
48  data Scalar = S !Double !Double deriving (Show,Eq)
```

The first number is the lower bound, the second the upper bound. The following function returns the upper bound:

```
49  upper (S _ u) = u
```

We bound operations on real numbers by first lifting them to operations on `Scalar` values and then bound the resulting operations by *enlarging* the bound to safe lower and upper bounds. An operation is *exact* if it does not involve any rounding (in which case there is no need to enlarge a bound).

The function that takes a `Scalar` with lower bound $l$ and upper bound $u$, then finds a safe lower bound on $l$ and a safe upper bound on $u$ is implemented as follows:

```
50  enlarge (S l u) = S (safeLowerBound l) (safeUpperBound u)
```

For convenience we provide a function to convert a number $x$ to a `Scalar` with $x$ as both lower and upper bound:

```
51  toScalar x = S x x
```

7.3. **Arithmetic on scalars.** We make `Scalar` an instance of the `Num` type class so that we can perform arithmetic on scalars (addition (`+`), subtraction (`-`), negation, multiplication (`*`) and non-negative integer exponentiation (`^`)).

```
52  instance Num Scalar where
```

If $x \in [l, u]$ then $-x \in [-u, -l]$; negation is exact on safe number so we do not need to enlarge the bound:

```
53      negate (S l u) = S (-u) (-l)
```

If $x \in [l, u]$ then $|x| \in [a, b]$, where $a = \max\{0, l, -u\}$ and $b = -\min\{0, l, -u\}$ (it is easy to check that this is correct regardless of the signs of $l$ and $r$). All operations involved are exact on safe numbers so we do not need to enlarge the bound:

```
54      abs (S l u) = S (maximum xs) (-minimum xs)
            where xs = [0, l, -u]
```

If $x \in [l, u]$ and $y \in [l', u']$, then $x + y \in [l + l', u + u']$. This operation is not exact so we enlarge the bound:

```
56      (S l u) + (S l' u') = enlarge (S (l + l') (u + u'))
```

If $x \in [l, u]$ and $y \in [l', u']$, then $x * y \in [a, b]$ where $a$ is the minimum of the numbers $\{l*l', l*u', u*l', u*u'\}$ and $b$ is the maximum of the same numbers. This operation is not exact so we enlarge the bound:

```
57      (S l u) * (S l' u') = enlarge (S (minimum xs) (maximum xs))
            where xs = [l*l', l*u', u*l', u*u']
```

The last two methods are required to complete the implementation of the `Num` instance (`fromInteger` provides implicit conversion of integer literals to `Scalar` values):

```
59      fromInteger   = toScalar . fromInteger
        signum (S l u) = error "S.signum: not defined"
```



In order to be able to divide `Scalar` values using (`/`) we must also add `Scalar` to the `Fractional` type class.

```
61  instance Fractional Scalar where
```

If $x \in [l, u]$ and if $l, u$ have the same sign, then the reciprocal is well defined for $x$ and $1/x \in [1/u, 1/l]$. This operation is not exact so we enlarge the bound:

```
62      recip (S l u) | l*u > 0   = enlarge (S (1/u) (1/l))
                     | otherwise = error "S.recip: not well-defined"
```

The last method is required; it provides implicit conversion of decimal literals to `Scalar` values:

```
64      fromRational = toScalar . fromRational
```

7.4. **Ordering of scalars.** In order to be able to compare `Scalar` values, e.g. using (`<`), we add `Scalar` to the `Ord` type class. If two bounds overlap we declare them incomparable and halt the program, otherwise comparison is implemented in the obvious way.

```
65  instance Ord Scalar where
        compare (S l u) (S l' u')
                | u < l'              = LT
                | l > u'              = GT
                | l == l' && u == u'  = EQ
                | otherwise           = error "S.compare: uncomparable"
```

## 8. Computation with polynomials

We show how to lift operations on polynomials (of degree $d-1$) to rectangle sets in $\mathbb{R}^d$ and then how to bound these operations.

8.1. **Representation of polynomials.** Polynomials are represented as a list of `Scalar` values (with the first element representing the constant coefficient). Hence what we refer to as a 'polynomial' of degree $d-1$ is actually a rectangle set in $\mathbb{R}^d$. In this section we lift operations on actual polynomials to such rectangles. We do not need to find any bounds on these lifts since this was already done implicitly in the previous section.

8.2. **Arithmetic with polynomials.** Add polynomials to the `Num` type class so that we can perform arithmetic operations on polynomials. (This implementation is a bit more general since it adds `[a]` to the `Num` type class for any type `a` in the `Num` type class.)

```
71  instance (Num a) => Num [a] where
```

Addition: $[c_1 + zq_1(z)] + [c_2 + zq_2(z)] = [c_1 + c_2] + z[q_1(z) + q_2(z)]$.

```
72      (c1:q1) + (c2:q2) = c1 + c2 : q1 + q2
        []      + p2      = p2
        p1      + []      = p1
```

Multiplication: $[c_1 + zq_1(z)] \cdot [c_2 + zq_2(z)] = [c_1 \cdot c_2] + z[c_1 \cdot q_2(z) + q_1(z) \cdot p_2(z)]$, where $p_2(z) = c_2 + zq_2(z)$.

```
75      (c1:q1) * p2@(c2:q2) = c1*c2 : [c1]*q2 + q1*p2
        _       * _          = []
```



The remaining methods are straightforward:

```
77    negate p      = map negate p
      fromInteger c = [fromInteger c]
      abs           = error "abs not implemented for polynomials"
      signum        = error "signum not defined for polynomials"
```

8.3. **Polynomial evaluation.** Evaluation of the polynomial $c + zq(z)$ at the point $t$ is done using the obvious recursion:

```
81    peval (c:q) t = c + t * peval q t
      peval []     _ = 0
```

8.4. **Norm of polynomial.** We use the $\ell^1$–norm on polynomials, i.e. $\|a_0 + \cdots + a_n z^n\| = |a_0| + \cdots + |a_n|$:

```
83    pnorm = sum . map abs
```

8.5. **Derivative of polynomial.** The derivative of $c + zq(z)$ is implemented using the recursion $D(c + zq(z)) = q(z) + zDq(z)$:

```
84    pderiv (c:q) = q + (0 : pderiv q)
      pderiv []    = []
```

## 9. Computation with analytic functions

We show how to lift operations on analytic functions to rectangle subsets in $X$ and how to bound these operations.

9.1. **The Function data type.** Functions in $X$ are represented as
$$f(z) = p(z) + z^d h(z),$$
where $p$ is a polynomial (not necessarily of degree less than $d$) and $\|h\| < K$, where $h \in X$. We refer to $p$ as the *polynomial part* of $f$ and $h$ is called the *error* of $f$. The value for the degree $d$ is specified in Appendix C.

The `Function` data type represents an analytic function on the above form (the first parameter is the polynomial part, the second parameter is the bound on the error):

```
86    data Function = F ![Scalar] !Scalar deriving (Show,Eq)
```

That is, `Function` represents rectangle subsets of $X$ of the form
$$\{a_0 + \cdots + a_n z^n + z^d h(z) \mid a_k \in A_k, k = 0, \ldots, n, \|h\| \in I\},$$
where $\{A_k\}$ and $I$ are intervals. Only the upper bound on the error term is needed so we do not take care to ensure that the lower bound is correct. Hence, the lower bound will be meaningless in general.

Note that we do allow $n \geq d$ in the above representation but in general we adjust our computations to ensure $n < d$. We call this operation *splitting*: if
$$f(z) = a_0 + \cdots + a_n z^n + z^d h(z),$$
with $n \geq k \geq d$, then we can split $f$ at degree $k$ into
$$f(z) = a_0 + \cdots + a_{k-1} z^{k-1} + z^d[a_k z^{k-d} + \cdots + a_n z^{n-d} + h(z)]$$
$$= p'(z) + z^d[r(z) + h(z)].$$



Thus the polynomial part of $f$ after splitting is $p'$ and the error is bounded by $\|r\| + \|h\|$ (by the triangle inequality). The implementation of this operation is:

```
87  split k (F p e) = let (p',r) = splitAt k p in F p' (e + pnorm r)
```

We will now lift operations on analytic functions to the above type of rectangles and then find bounds on these operations.

9.2. **Arithmetic with analytic functions.** In what follows we let $f_i(z) = p_i(z) + z^d h_i(z)$ for $i = 1, 2, 3$, and let $f_1 \diamond f_2 = f_3$ where $\diamond$ is the operation under consideration.

Make `Function` an instance of the `Num` type class so that we can perform arithmetic operations on functions (addition (+), subtraction (-), negation, multiplication (*) and non-negative integer exponentiation (^)).

```
88  instance Num Function where
```

Addition of two functions is performed by adding the polynomial part and the error separately, $p_3 = p_1 + p_2$ and $h_3 = h_1 + h_2$, so that $\|h_3\| \leq \|h_1\| + \|h_2\|$ by the triangle inequality:

```
89      (F p1 e1) + (F p2 e2) = F (p1 + p2) (e1 + e2)
```

Multiplication of two analytic functions is given by the equation

$$f_1(z)f_2(z) = p_1(z)p_2(z) + z^d \left[ p_1(z)h_2(z) + p_2(z)h_1(z) + z^d h_1(z) h_2(z) \right],$$

so that $\|h_3\| \leq \|p_1\|\|h_2\| + \|p_2\|\|h_1\| + \|h_1\|\|h_2\|$. To ensure that the degree of the polynomial part does not increase too much we split it at degree $d+1$:[17]

```
90      (F p1 e1) * (F p2 e2) = split (d+1) (F (p1*p2) e3)
            where e3 = e2*pnorm p1 + e1*pnorm p2 + e1*e2
```

The negation of $f$ is $-p(z) + z^d(-h(z))$ but the error is unchanged since we only keep a bound on its norm:

```
92      negate (F p e) = F (negate p) e
```

The remaining methods are required to complete the implementation of the `Num` instance (`fromInteger` provides implicit conversion of integer numerals to `Function` values):

```
93      fromInteger c = F [fromInteger c] 0
        abs           = error "abs not implemented for Function"
        signum        = error "signum not defined for Function"
```

9.3. **Norm of analytic functions.** The triangle inequality gives $\|p(z) + z^d h(z)\| \leq \|p\| + \|h\|$ (since $|z| < 1$):

```
96  norm (F p e) = pnorm p + e
```

The norm on the Cartesian product $Y = X \times X$ is $\|(f, g)\| = \|f\| + \|g\|$:

```
97  sumnorm (f,g) = norm f + norm g
```

The distance induced by the norm on $Y$:

```
98  dist (f,g) (f',g') = sumnorm (f-f',g-g')
```

---

[17] We choose to split at degree $d+1$ (instead of the perhaps more natural choice of degree $d$) because the division by $z$ in the definition of the operator $T$ would otherwise cause $g_8$ to have degree at most $d-2$.



9.4. **Differentiation.** The implementation of differentiation of $f \in X$ is complicated by the use of the $\ell^1$–norm on $X$, since $\|f\| < \infty$ does not imply that $\|Df\| < \infty$. This problem is overcome by only computing the derivative of functions restricted to a disk of radius strictly smaller than one. That is, we need to know a-priori that the function we are differentiating only will be evaluated on this smaller disk. Usually we get this information from the fact that we compute derivatives like $Df_1 \circ f_2$ and we have bounds on the image of $f_2$.

Given $f \in X$ we will estimate $Df\,|_{\{|z|<\mu\}}$ where $\mu < 1$. If $f(z) = p(z) + z^d h(z)$, then $Df(z) = Dp(z) + dh(z)z^{d-1} + z^d Dh(z) = p_1(z) + z^d h_1(z)$. Here we are faced with the problem that we only know the norm of $h$ so all we can say about the polynomial part is that $p_1(z) = Dp(z) + sdz^{d-1}$, where $s \in [-\|h\|, \|h\|]$.

Let $h(z) = \sum a_k z^k$, then the error can be crudely approximated as follows:

$$\|Dh(z)\,|_{\{|z|<\mu\}}\| = \|Dh(\mu z)\| = \sum_{k \geq 1} k\mu^{k-1}|a_k| \leq \sup_{k \geq 1} k\mu^{k-1} \|h\|$$

$$\leq \|h\| \sum_{k \geq 1} k\mu^{k-1} = \frac{\|h\|}{(1-\mu)^2}$$

Putting all this together we arrive at the following implementation:

```
99   deriv mu (F p e) | mu < 1    = F p1 e1
                     | otherwise = error "deriv: mu is not < 1"
        where p1 = pderiv p + [S (-s) s] * [0,1]^(d-1)
              e1 = e / (1 - mu)^2
              s  = fromIntegral d * upper e
```

Note that $\mu$ is passed as a parameter by the caller of this function (it is not a constant). As mentioned earlier, usually this function is used to compute expressions like $Df_1 \circ f_2$ in which case $\mu$ will be an upper bound on the radius of a disk containing the image of $f_2$.

9.5. **Composition.** The implementation of composition of analytic functions $f_1 \circ f_2$ is split up into two parts. First we consider the special case when $f_1 = p_1$ is a polynomial, then we treat the general case.

Polynomials are defined on all of $\mathbb{C}$ so the composition $p_1 \circ f_2$ is always defined. If $p_1(z) = c + zq(z)$ then we may use the recursion suggested by $p_1 \circ f_2(z) = c + f_2(z) \cdot q \circ f_2(z)$:

```
104  compose' (c:q) f2 = (F [c] 0) + f2 * compose' q f2
```

The recursion ends when the polynomial is the zero polynomial, in which case $p_1 \circ f_2 = 0$:

```
105  compose' []     _  = 0
```

In the general case we have to take care to ensure that the image of $f_2$ is contained in the domain of $f_1$ for the composition to be well-defined. A sufficient condition for this to hold is $\|f_2\| < 1$ since the domain of $f_1$ is the unit disk. Under this assumption we compute $f_1 \circ f_2 = p_1 \circ f_2 + (f_2)^d \cdot h_1 \circ f_2$. These two terms are split at degree $d+1$ to get $p_1 \circ f_2(z) = \tilde{p}_1(z) + z^d \tilde{h}_1(z)$ and $f_2(z)^d = \tilde{p}_2(z) + z^d \tilde{h}_2(z)$.[18] Then $f_1 \circ f_2(z) = p_3(z) + z^d h_3(z)$ with $p_3 = \tilde{p}_1 + \tilde{p}_2 \cdot h_1 \circ f_2$ and $h_3 = \tilde{h}_1 + \tilde{h}_2 \cdot h_1 \circ f_2$.

---

[18] See the footnote near the definition of multiplication of analytic functions for an explanation of the choice of degree $d+1$.



Only the norm of $h_1$ is given so from this we can only draw the conclusion that $p_3(z) = \tilde{p}_1(z) + s \cdot \tilde{p}_2(z)$ for $s \in [-\|h_1\|, \|h_1\|]$ ($s$ is in fact a function but we may think of it as a constant since we are really computing with sets of polynomials). The error is approximated using the triangle inequality, $\|h_3\| \leq \|\tilde{h}_1\| + \|\tilde{h}_2\|\|h_1\|$.

```
106 compose (F p1 e1) f2 | norm f2 < 1 = F (p1' + [s]*p2') (e1' + e1*e2')
                        | otherwise   = error "compose: |f2| is too large"
       where (F p1' e1') = split (d+1) (compose' p1 f2)
             (F p2' e2') = split (d+1) (f2^d)
             s           = S (-upper e1) (upper e1)
```

The term $s \cdot \tilde{p}_2(z)$ can introduce devastating errors into the computation since $s$ lies in an interval which has positive upper bound and a negative lower bound (if $\tilde{p}_2$ has a coefficient with small error but a large magnitude relative to $s$, then after multiplying with $s$ that coefficient will have an error that is bigger than the magnitude of the coefficient). We work around this problem by choosing the degree $d$ large, since this tends to make the term $s$ smaller. Another way to deal with this problem is to include a "general error" term in the representation of analytic functions (see [10]).

9.6. **Derivative of the composition operator.** Let $S(f,g) = f \circ g$, then the derivative is given by
$$DS_{(f,g)}(\delta f, \delta g) = Df \circ g \cdot \delta g + \delta f \circ g.$$

Note that when computing $Df$ we must specify as a first parameter the radius of a disk strictly contained in the unit disk to which $Df$ is restricted (see Section 9.4). In the present situation we know that the image of $g$ is contained in a disk with radius $\|g\|$, so $Df$ only needs to be evaluated on the disk of radius $\|g\|$:

```
111 dcompose f g df dg = (deriv (norm g) f `compose` g) * dg + (df `compose` g)
```

9.7. **Division by $z$.** If $f(z) = a_1 z + \cdots + a_n z^n + z^d h(z)$, then
$$f(z)/z = a_1 + \cdots + a_n z^{n-1} + z^{d-1} h(0) + z^d \tilde{h}(z),$$

where $|h(0)| \leq \|h\|$ and $\|\tilde{h}\| \leq \|h\|$. Since we do not know the value of $h(0)$ we estimate the coefficient it with $s \in [-\|h\|, \|h\|]$. We think of this operation as a "left shift", whence the name of this function:

```
112 lshift (F (c:q) e) = F (q + [0,1]^(d-1) * [S (-upper e) (upper e)]) e
    lshift (F [] e)    = F ([0,1]^(d-1) * [S (-upper e) (upper e)]) e
```

If the polynomial part of $f$ has a constant coefficient $a_0 \neq 0$ then this function will not return the correct result, so we take care to only use it when we know that $a_0 = 0$.

9.8. **Point evaluation.** If $f(z) = p(z) + z^d h(z)$, then $f(t) = p(t) + t^d \cdot s$ for some $s \in [-\|h\|, \|h\|]$. We also check that $t$ is in the unit disk otherwise the program is terminated with an error:

```
114 eval (F p e) t | abs t < 1 = peval p t + t^d * (S (-upper e) (upper e))
                   | otherwise = error ("eval: not in domain t=" ++ show t)
```

Note that the further away $t$ is from 0, the more error is introduced in the evaluation. For $t = 0$ the error term has no influence on the evaluation.



9.9. **Scaling.** As a convenience we define operators to scale an analytic function by a scalar on the on the right. The precedence for these operators are the same as for their 'normal' counterparts.

Multiplication satisfies $[p(z) + z^d h(z)] \cdot x = x \cdot p(z) + z^d[x \cdot h(z)]$ and division is handled similarly. Note that the error term is affected:

```
infixl 7 .*, ./
(F p e) .* x = F (p * [x]) (e * abs x)
(F p e) ./ x = F (p * [1/x]) (e / abs x)
```
116

9.10. **Polynomial approximation.** Let $f(z) = p(z) + z^d h(z)$. To approximate $f$ by a polynomial we first discard the error term $z^d h(z)$, then we disregard the errors in the coefficients of $p$. That is, for $p(z) = a_0 + \cdots + a_n z^n$ with $a_k \in [a_k^-, a_k^+]$ we replace $a_k$ with the mean $\tilde{a}_k = (a_k^- + a_k^+)/2$ (we 'collapse' the bounds on $a_k$). Finally, we lift this operation to pairs of functions:

```
approx (f,g) = (approx' f, approx' g)
    where approx' (F p _)  = F (map (toScalar . collapse) p) 0
          collapse (S l u) = (l+u)/2
```
119

9.11. **Operator norm.** Let $\xi_k(z) = z^k$ so that $\{\xi_k\}_{k \geq 0}$ is a basis for $X$. A basis for $Y$ is $\{\eta_k\}_{k \geq 0}$, where $\eta_{2k} = (\xi_k, 0)$ and $\eta_{2k+1} = (0, \xi_k)$. This set is implemented as follows:

```
basis = interleave (zip basis' (repeat 0)) (zip (repeat 0) basis')
    where basis' = map xi [0..]
          xi k   = F (replicate k 0 ++ [1]) 0
```
122

**Proposition 9.1.** *If* $L : Y \to Y$ *is a linear and bounded operator, then*

$$\|L\| = \max\{\|L\eta_0\|, \ldots, \|L\eta_{2d-1}\|, \sup_{h \in B_d} \|L(h,0)\|, \sup_{h \in B_d} \|L(0,h)\|\},$$

*where* $B_d = \{z^d h(z) \mid \|h\| < 1\}$.

This is a consequence of using the $\ell^1$–norm on $X$.

Given a linear operator `op` acting on a list of tangent vectors[19], we estimate the operator norm by applying it to the first $2d$ basis vectors and to the sets $B_d \times 0$ and $0 \times B_d$. Then we compute the upper bound of the norm of the results and take the maximum:

```
opnorm op = maximum $ map (upper . sumnorm) $ op tangents
    where tangents = (F [] 1,0) : (0,F [] 1) : take (2*d) basis
```
125

Note that $B_d$ is represented by the set of functions with no polynomial part and an error bounded by 1, which is the same as `F [] 1`.

---

[19] The linear operator acts on a sequence of tangent vectors since this is how we have implemented the derivative of the main operator.



### 9.12. Construction of balls.
We cannot exactly represent arbitrary balls in $X$ with the `Function` type. Instead we construct a rectangle set which is guaranteed to contain the ball.

Thus, a bound on a ball of radius $r$ centered on an analytic function (in our case it is always a polynomial, i.e. `e=0`) can be implemented as follows:

```
127  ball r (f,g) = (ball' r f, ball' r g)
       where ball' r (F p e) = F (map (+ S (-r) r) p) (e + toScalar r)
```

### 9.13. Newton's method.
This is our variant of Newton's method on $Y$:
$$(f,g) \mapsto (M - I)^{-1}(M - T)(f,g),$$
where $M$ is a $2d \times 2d$ matrix passed as the first parameter. The second parameter is $T(f,g)$ and the third parameter is $(f,g)$. When lifting $M$ into $Y$ we project the error term to zero by letting $s = 0$ and when lifting $(M - I)^{-1}$ we preserve the error term by letting $s = 1$ (see below for how this lifting is done):

```
129  newton m (tf,tg) fg = fg'
       where (mf,mg) = liftPolyOp 0 (apply m) fg
             fg'     = liftPolyOp 1 (solve $ subtractDiag m 1) (mf-tf,mg-tg)
```

Let $f = p(z) + z^d h(z)$ where $\deg p < d$ and let $A$ be the linear operator represented (in the basis $\{z^k\}$) by the infinite matrix
$$\begin{pmatrix} M & 0 \\ 0 & sI \end{pmatrix}$$
where $M$ is a $d \times d$ matrix and $I$ is the infinite identity matrix. We lift the linear operator $A$ into $X$ by
$$Af(z) = Mp(z) + z^d(s \cdot h(z)).$$

The following function implements this lifting into $Y$. We split $f$ and $g$ to ensure that their degrees are at most $d - 1$, and since our linear algebra routines require its input in one vector we interleave the polynomial parts. Also, instead of passing $M$ we pass a linear operator `op` which allows us to use one function to lift matrix multiplication (`apply`) and solution of linear equations (`solve`):

```
132  liftPolyOp s op (f,g) = (F pf' (s*ef), F pg' (s*eg))
       where fg@(F pf ef, F pg eg) = (split d f, split d g)
             (pf',pg') = uninterleave $ op (interleavePoly fg)
```

When interleaving the polynomial parts of two functions we first pad the polynomials with zeros to ensure their lengths are exactly $d$ (e.g. $a_0 + a_1 z$ is padded to $a_0 + a_1 z + 0z^2 + \cdots + 0z^{d-1}$). Hence the resulting vector always has length $2d$:

```
135  interleavePoly (F p _, F q _) = interleave (pad p) (pad q)
       where pad x = take d $ x ++ (repeat 0)
```

### Appendix A. Linear algebra routines

In this section we implement a simple linear algebra library to compute matrix-vector products and to solve linear equations.

A matrix is represented as a list of its rows and a row is a list of its elements. A vector is just a list of elements (we think of them as column vectors). This is a very simplistic library so no checking is done to ensure that matrices have the



correct dimensions (e.g. it is quite possible to create a 'matrix' with rows of differing lengths).

A.1. **Matrix-vector product.** Computing the matrix-vector product $Mx$ is fairly straightforward:

```
137  apply m x = map (dotProduct x) m
```

The dot product of vectors $a$ and $b$:

```
138  dotProduct a b = sum $ zipWith (*) a b
```

Note that if `a` and `b` have different lengths, then the above function will treat the longer vector as if it had the same length as the shorter.

A.2. **Linear equation solver.** The following function solves the linear system of equations $Mx = b$. It is a simple wrapper around a function which solves a linear system given an augmented matrix.

```
139  solve m b = solveAugmented $ augmentedMatrix m b
```

The augmented matrix for $M$ and $b$ is $\begin{bmatrix} M & b \end{bmatrix}$, i.e. the matrix with $b$ appended as the last column of $M$:

```
140  augmentedMatrix = zipWith (\x y -> x ++ [y])
```

We now implement the linear equation solver which takes an augmented matrix as its only parameter. It is implemented using Gaussian elimination with partial pivoting. The only novelty compared with a traditional imperative implementation is that we solve the equations recursively.

Given a $n \times (n+1)$ augmented matrix $M'$ first perform partial pivoting, i.e. the row whose first element has the largest magnitude is moved to the top to form the matrix $M$. Assuming that we already have the solution for $x_2, \ldots, x_n$ we can compute $x_1 = (m_{1,n+1} - \sum_{j=2}^{n} m_{1j} x_j)/m_{11}$ and we are done. The solution for $x_2, \ldots, x_n$ is found recursively as follows: perform a Gaussian elimination on $M$ to ensure that all rows except the first start with a zero to get a matrix $\tilde{M}$. Throw away the first row and column of $\tilde{M}$ to get a $(n-1) \times n$ matrix $N'$ and solve the linear system with augmented matrix $N'$. The solution to this system is $x_2, \ldots, x_n$.

```
141  solveAugmented [] = []
     solveAugmented m' = (last m1t - dotProduct m1t x) / m11 : x
         where m@((m11:m1t):_) = partialPivot m'
               x = solveAugmented $ eliminate m
```

Partial pivoting is done by first finding a list of all possible ways to split the matrix $M$ into a top and a bottom half. This list is searched for the split which has a maximal first element in the bottom half. The maximal split is then reassembled into one matrix by moving the top row of the bottom half to the top of the matrix.

```
145  partialPivot m = piv:mtop ++ mbot
         where (mtop,piv:mbot) = maximumBy comparePivotElt (splits m)
               comparePivotElt (_,(a:_):_) (_,(b:_):_) = compare (abs a) (abs b)
```

The following routine uses Gaussian elimination to ensure that the first element of all rows except the first starts with a zero. That is, we add a suitable multiple of the first row to the other rows one at a time:



```
148  eliminate ((m11:m1t):mbot) = foldl appendScaledRow [] mbot
         where appendScaledRow a (r:rs) = a ++ [scaleAndAdd (-r/m11) m1t rs]
               scaleAndAdd s a b = zipWith (+) (map (*s) a) b
```

*Remark* A.1. When using the above linear equation solver with matrices and vectors over intervals (of type `Scalar`) there is a question of what the 'solution' represents. As always, we are computing *bounds* on solutions: if $M$ is in some rectangle set $[M]$ of matrices and $b$ is in some rectangle set $[b]$ of vectors, then the above routine will compute a rectangle set $[x]$ such that if $x$ is a solution to $Mx = b$ then $x \in [x]$.

Note that our solver will compute rather loose bounds on the solution set, see e.g. [8] for ways of finding sharper bounds.

## Appendix B. Supporting functions

Given a square matrix $M$ and a number $x$ compute $M - xI$, i.e. subtract $x$ from every diagonal element of $M$:

```
151  subtractDiag m x = foldl f [] (zip m [0..])
         where f m' (r,k) = let (h,t:ts) = splitAt k r
                            in m' ++ [h ++ [t-x] ++ ts]
```

Given a list, return all possible ways to split the list in two:

```
154  splits x = splits' [] x
         where splits' _ []        = []
               splits' x y@(yh:yt) = (x,y) : splits' (x ++ [yh]) yt
```

Interleave two lists `a` and `b`, i.e. construct a new list by taking the first element from `a`, then the first element from `b` and repeating for the remaining elements.

```
157  interleave a b = concat $ zipWith (\x y -> [x,y]) a b
```

Perform the 'inverse' of the above function, i.e. take a list `c` and construct a pair of lists `(a,b)` such that `interleave a b = c`:

```
158  uninterleave = unzip . pairs
```

Given a list, partition it into pairs of adjacent elements:

```
159  pairs [] = []
     pairs (x:y:rest) = (x,y) : pairs rest
     pairs _ = error "list must have even length"
```

## Appendix C. Input to the main program

The degree of the error term in our representation of analytic functions:

```
162  d = 13 :: Int
```

The radius of the ball on which $\Phi$ is a contraction:

```
163  beta = 1.0e-7 :: Double
```

The radii for the domains of $\phi$ and $\psi$:

```
164  sf = 2.2 :: Scalar
     sg = 0.5 :: Scalar
```

The initial guess for the fixed point:

```
166  guess = (F [-0.75, -2.5] 0, F [6.2,-2.1] 0)
```



Appendix D. Running the main program

This document contains all the Haskell source code needed to compile the program into an executable. Given a copy of the LaTeX source of this document (assuming the file is named `lmca.lhs`), use the following command to compile it:[20]

```
ghc --make -O2 lmca.lhs
```

This produces an executable called `lmca` (or `lmca.exe` if you are using Windows) which when called will execute the `main` function.

Here is the output of running the main program:

```
radius     = 1.0e-7,
|Phi(f)-f| < 4.830032057738462e-9,
|DPhi|     < 0.1584543808202988
```

This output was taken from a sample run using GHC 6.12.1 on Mac OS X 10.6.2. The running time on an 1.8 GHz Intel Core 2 Duo was less than 10 seconds.

Appendix E. Haskell mini-reference

This section introduces some of the features and syntax of Haskell to help anybody unfamiliar with the language read the source code. It is assumed that the reader has some prior experience with an imperative language (Java, C, etc.) but is new to functional programming languages. Table 1 below collects examples of Haskell syntax used in the source code and can be used to look up unfamiliar expressions. For more information on the Haskell language go to http://haskell.org.

Haskell is a *functional* language. Such languages differ from imperative languages in several significant ways, e.g.: there are no control structures such as `for` loops and data is immutable so there is no concept of variables (memory locations) that can be written to.

Basic types include: *booleans* (`True`, `False`), *numbers* (e.g. `-1`, `2.3e3`, integers of any magnitude are supported), *tuples* (e.g. `(1,'a',0.3)`, elements can have different types), and *lists* (e.g. `[1,2,3]`, all elements must have the same type). Functions are on the same level as basic types so they can e.g. be passed as parameters to other functions.

Functions are defined like `f` *parameters* = *expression* where `f` is the function name and there can be zero or more parameters. Note that there are no parentheses around parameters and that parameters are separated by spaces. Function calls have very high precedence, so `f x^2` is the same as `(f x)^2`, *not* `f (x^2)`. The keywords `let .. in` and `where` can be used to bind expressions to function-local definitions (i.e. local functions or variables).

New data types can be defined using the `data` construct. For example, `data Interval = I Double Double` defines a type called `Interval` which consists of two double precision floating point numbers (i.e. the endpoints of the interval). New values of this type are constructed using the *value constructor* which we called `I`, e.g. `I 0 1` defines the unit interval.

Functions can be defined with *pattern matching* on built-in and custom data types. For example, `len (I a b) = abs (b-a)` defines a function `len` which returns the length of an interval (for the custom data type `Interval`).

---

[20] The GHC compiler can be downloaded for free from http://haskell.org.



We often use pattern matching on lists, where `[]` matches the empty list and `(x:xs)` matches a list with a least one element and binds the first element to `x` and the rest to `xs` (read as plural of `x`). The notation `v@(x:xs)` can be used to bind the entire list to `v` on a match.

The notation `_` may be used to match anything without binding the match to a variable, e.g. `firstZero (x:_) = x == 0` defines a function which returns `True` if the first element of a non-empty list is equal to zero (and throws an exception if called on the empty list `[]`).

*Type classes* are a way of declaring that a custom data type supports a certain predefined collection of functions and also allow for 'overloading' of functions (and operators, which can be turned into functions as noted in the example for `(+)` in Table 1). We only mention type classes because we come across them when implementing `Scalar` and `Function`. The pre-defined type classes we use are `Num` (for `(+)`, `(-)`, `(*)`, `(^)`, `abs`), `Fractional` (for `(/)`, `(^^)`), `Eq` (test for equality), and `Show` (for conversion to strings).

TABLE 1. Examples of Haskell syntax used in the source code.

| Expression | Description |
|---|---|
| `f1 x = 2*x` | define a function `f1` which doubles its argument |
| `f2 x y = x+y` | define a function `f2` which adds its two arguments |
| `f1 3` | apply `f1` to `3` (=6) |
| `f2 3 4` | apply `f2` to `3` and `4` (= 7) |
| `f2 2 (f1 3)` | apply `f2` to `2` and `6` (the result of `f1 3`) (= 8) |
| `f2 2 $ f1 3` | same as above (the operator `$` is often used in this way to avoid overuse of parentheses) |
| `f2 2 f1 3` | error (this means: compute `f2 2 f1` and apply the result to `3`, but `2+f1` does not make sense) |
| `\x -> 2*x` | define the *anonymous function* $x \mapsto 2x$ |
| `f2 3` | apply `f2` to `3` (= the function `\x -> 3+x`) |
| `f1 . f2 3` | function composition (= the function `\x -> 2*(3+x)`) |
| `` 3 `f2` 4 `` | turn function (in backticks) into an operator (=7) |
| `(+) 3 4` | turn operator (in parentheses) into a function (=7) |
| `(3*)` | fix first parameter to `3` (= the function `\x -> 3*x`) |
| `g x | x<0  = -1`<br>`    | x>0  =  1`<br>`    | x==0 =  0` | define the sign function `g` using *guards* (the `|` symbols) |
| `num [] = 0`<br>`num (_:xs)`<br>`    = 1+num xs` | define a function which counts the number of elements in a list using pattern matching |
| `[1,2]` | a list (all elements must have the same type) |
| `[]` | the empty list |
| `[1..]` | the list of all positive integers |
| `[2..5]` | list enumeration with bounds (=`[2,3,4,5]`) |
| `1 : [2,3]` | append element to beginning of list (=`[1,2,3]`) |

*Continued on the next page...*



| Expression | Description |
|---|---|
| `[1,2,3] !! 0` | access list elements by zero-based index ($=1$) |
| `[1,2] ++ [3]` | concatenate two lists ($=$`[1,2,3]`) |
| `'a'` | a character |
| `"abc"` | a string, i.e. a list of characters ($=$`['a','b','c']`) |
| `2^3` | non-negative integer exponentiation ($=8$) |
| `2^^(-1)` | integer exponentiation ($=0.5$) |
| `('a',2)` | a pair (the elements need not have the same type) |
| `fst ('a',2)` | access first element in a pair ($=$`'a'`) |
| `snd ('a',2)` | access second element in a pair ($=2$) |
| `map f1 [1..]` | apply `f1` to all elements in the list ($=$`[2,4,8,..]`) |
| `[f2 a b \| a <- [1,2], b <- [3..5]]` | list comprehension, i.e. $\{f_2(a,b) \mid a \in \{1,2\}, b \in \{3,4,5\}\}$ ($=$`[4,5,6,5,6,7]`) |
| `foldl f2 1 [3,5]` | fold left over list (compute `f2 1 3` $=4$, then `f2 4 5`) ($=9$) |
| `iterate f1 1` | compute orbit of `1` under `f1` ($=$`[1,2,4,8,..]`) |
| `maximum [1,4,2]` | return maximum element in a list ($=4$) |
| `maximumBy f x` | as above, but using `f` to compare elements of the list `x` |
| `minimum [1,4,2]` | return minimum element in a list ($=1$) |
| `splitAt 2 [1,4,2]` | split list in two at given index ($=$`([1,4],[2])`) |
| `take 3 [7..]` | take the first 3 elements from the list ($=$`[7,8,9]`) |
| `zip [1..] [3,4]` | join two lists into a list of pairs ($=$`[(1,3),(2,4)]`) |
| `unzip [(1,3),(2,4)]` | 'inverse' of `zip` ($=$`([1,2],[3,4])`) |
| `zipWith f2 [1..] [3,4]` | like `zip`, but use `f2` to join elements ($=$`[4,6]`) |
| `repeat 0` | infinite list with one element repeated ($=$`[0,0,0,..]`) |
| `replicate 3 0` | finite list with one element repeated ($=$`[0,0,0]`) |
| `sum [3,-1,4]` | sum elements in list ($=6$) |
| `transpose m` | the transpose of the matrix `m` (`m` is a list of lists) |
| `putStrLn "hi"` | print `hi` to standard out and and append a new line |
| `show 1.2` | turn the number `1.2` into the string `"1.2"` |
| `error "ohno"` | abort program with error message `ohno` |

*E-mail address*: `winckler@kth.se`